\documentclass[11pt]{amsart}
\usepackage{amssymb}
\usepackage{amsfonts}
\usepackage{graphicx}
\usepackage{color}

\newtheorem{theorem}{Theorem}[section]

\newtheorem{lemma}[theorem]{Lemma}
\newtheorem{corollary}[theorem]{Corollary}

\renewenvironment{proof}{{\noindent\bf Proof.}}{\hfill $\Box$\par\vskip3mm}

\begin{document}

\def\reals{{\mathbb R}}
 \def\ch{{\mathcal H}}
 \def\cA{{\mathcal A}}
 \def\cD{{\mathcal D}}
 \def\cK{{\mathcal K}}
 \def\cC{{\mathcal C}}
 \def\cN{{\mathcal N}}
 \def\cR{{\mathcal R}}
 \def\cS{{\mathcal S}}
 \def\cT{{\mathcal T}}
 \def\cV{{\mathcal V}}
 \def\ta{{\mathcal T}_{\subset}}
 \def\cI{{\mathcal I}}
 \def\bC{{\bf C}}
 \def\axis{{\bf A}}
 \def\fibr{{\bf H}}
 \def\ba{{\bf a}}
 \def\bb{{\bf b}}
 \def\bc{{\bf c}}
 \def\be{{\bf e}}
 \def\d{{\delta}} 
 \def\ci{{\circ}} 
 \def\e{{\epsilon}} 
 \def\l{{\lambda}} 
 \def\L{{\Lambda}} 
 \def\m{{\mu}} 
 \def\n{{\nu}} 
 \def\o{{\omega}} 
 \def\s{{\sigma}} 
 \def\v{{\varphi}} 
 \def\a{{\alpha}} 
 \def\b{{\beta}} 
 \def\p{{\partial}} 
 \def\r{{\rho}} 
 \def\ra{{\rightarrow}} 
 \def\lra{{\longrightarrow}} 
 \def\g{{\gamma}} 
 \def\D{{\Delta}} 
 \def\La{{\Leftarrow}} 
 \def\Ra{{\Rightarrow}} 
 \def\x{{\xi}} 
 \def\c{{\mathbb C}} 
 \def\z{{\mathbb Z}} 
 \def\2{{\mathbb Z_2}} 
 \def\q{{\mathbb Q}} 
 \def\t{{\tau}} 
 \def\u{{\upsilon}} 
 \def\th{{\theta}} 
 \def\la{{\leftarrow}} 
 \def\lla{{\longleftarrow}} 
 \def\da{{\downarrow}} 
 \def\ua{{\uparrow}} 
 \def\nwa{{\nwtarrow}} 
 \def\swa{{\swarrow}} 
 \def\nea{{\netarrow}} 
 \def\sea{{\searrow}} 
 \def\hla{{\hookleftarrow}} 
 \def\hra{{\hookrightarrow}} 
 \def\sl{{SL(2,\mathbb C)}} 
 \def\ps{{PSL(2,\mathbb C)}} 
 \def\qed{{\hfill$\diamondsuit$}} 
 \def\pf{{\noindent{\bf Proof.\hspace{2mm}}}} 
 \def\ni{{\noindent}} 
 \def\sm{{{\mbox{\tiny M}}}} 
 \def\sc{{{\mbox{\tiny C}}}} 

\title{Climbing a Legendrian mountain range without Stabilization}

\begin{abstract}
We introduce a new braid-theoretic framework with which to understand the Legendrian and transversal classification of knots, namely a {\em Legendrian Markov Theorem without Stabilization} which induces an associated {\em transversal Markov Theorem without Stabilization}.  We establish the existence of a nontrivial knot-type specific Legendrian and transversal MTWS by enhancing the Legendrian mountain range for the $(2,3)$-cable of a $(2,3)$-torus knot provided by Etnyre and Honda, and showing that elementary negative flypes allow us to move toward maximal $tb$ value without having to use Legendrian stabilization.  In doing so we obtain new ways to visualize convex tori and Legendrian divides and rulings, using tilings and braided rectangular diagrams.
\end{abstract}

\author{Douglas J. LaFountain and William W. Menasco}
\thanks{2000 \textit{Mathematics Subject Classification}. Primary 57M25, 57R17;
Secondary 57M50}
\keywords{Contact structures, convex surfaces, braids, Legendrian, transverse}
\date{\today}

\maketitle
\section{Introduction}
For closed braid representations of topological links in $S^3$ the {\em Markov Theorem
without Stabilization} (MTWS) \cite{[BM4]} states that for a fixed braid index $n^b$ there are a finite
number of ``modeled'' isotopies (dependent only on $n^b$ and not on link type) that take any
oriented link represented as an $n^b$-braid to a representative of minimal braid index
without the need for increasing the braid index via stabilization---an isotopy essential
in the classical Markov Theorem for closed braid equivalence.  Moreover, once at minimal
braid index the MTWS states that there are again a finite number of ``modeled'' isotopies
(again, dependent only on the value of the braid index) that allow us to jump between
conjugacy classes of minimal index.  These isotopies, which will grow in number as $n^b$ grows,
make up the MTWS {\em calculus} for closed braids \cite{[BM2]}.\\

Deciphering the inner structure of the MTWS calculus is a rich area for research.
To give an example, for $n^b = 3$ the MTWS calculus is made up of four closed
braid isotopies: positive braid preserving flypes; negative braid preserving flypes;
positive destabilizations; and, negative destabilizations \cite{[BM3],[BM6]}.
The MTWS calculus may also be calculated for specified link classes:  the calculus for the unlink
and torus knots utilizes only exchange moves and positive/negative destabilizations \cite{[BM1],[M1]}.
(For readers not familiar with these modeled braid isotopies, the introduction of \cite{[BM2]} is a suitable
source to consult.)\\

If we specialize the calculus to transversal $3$-braid knots in the standard-symmetric contact
structure for ${\mathbb R}^3 ( \subset S^3)$---the kernel of the $1$-form $dz + r^2 d \theta$---we
get an intriguing glimpse of the structure within the MTWS calculus.  Specifically, although positive
destabilizations, exchange moves, and positive braid preserving flypes are transverse isotopies,
closed braids which admit a negative braid preserving flype may not be transversally simple---classified
by their self-linking number \cite{[BW]}.  In fact, transversal knots having a closed braid representation
of minimal index $3$ which admit a negative flype but not a positive flype are the first explicit
examples of transversally non-simple knots \cite{[BM5]}.\\

The purpose of this note is to give evidence to a rich interplay between the structure of
the MTWS calculus and the classification structure of Legendrian and transversal
knot classes.  This evidence comes from the synthesis of three different lines of inquiry: characteristic
foliations of convex tori in a contact structure \cite{[E1],[G],[H]}; standard tiling of tori coming from
singular braid foliations \cite{[BM7],[BF]}; and the representation of Legendrian
knots by rectangular diagrams \cite{[MM]}.\\

This synthesis is brought to bear on the Etnyre-Honda ``Legendrian mountain range'' classification of the
Legendrian and transversal classes of the $(2,3)$ cabling of the $(2,3)$-torus knot---the
first implicit example of a non-simple transversal knot \cite{[EH]}.
Our initial result is an {\em enhanced mountain range classification} paradigm---a calculus structure
imposed over the mountain range that enables one to ``climb'' the mountain
without having to stabilize.  The implication of this enhanced classification structure is
that buried within the structure of the MTWS calculus lies a {\em Legendrian MTWS}.
Next, employing a result of Epstein, Fuchs and Meyer, 
this enhanced mountain range will ``collapse'' to an {\em enhanced mountain trail}
yielding the classification of transversal classes with a {\em transversal MTWS} calculus imposed---non-tranversal isotopies that allow us to jump between tranversal classes
without negative (non-transversal) stabilization.  By Bennequin's classical transversal result \cite{[B]}
and its Legendrian analogy \cite{[MM]}, the isotopies in these specialized Legendrian MTWS and transversal MTWS
can be represented by isotopies on rectangular closed braids and closed braids, respectively.
Finally, there remain modeled isotopies of the topological MTWS that are transversal
isotopies---including positive destablilization and exchange moves---and can be utilized to
move between closed braid representatives of the same transversal class without
increasing the braid index.\\

It is reasonable to conjecture that this stratification in the MTWS
calculus specialized to
the $(2,3)$ cabling of the $(2,3)$-torus knot is prototypical. That is,
coming out of the topological MTWS calculus
there is a Legendrian MTWS calculus for all Legendrian links
which will collapses to a transversal MTWS calculus; and what
remains will be the transversal isotopies of the MTWS. We advocate that
understanding this
stratifying structure is important to the study of contact geometry knot
theory and a fertile area
of research. In this vein, we take note of H. Matsuda's recent work on
Stalling's links and their connection
with additional structural aspects of the MTWS calculus landscape
\cite{[Ma]}.
From a rudimentary understanding of Matsuda's calculation for
determining the MTWS calculus for
$4$-braids or a review of such modeled isotopies as that in Figure 8 of
\cite{[BM4]}
we are lead to insert a cautionary remark. Not every Legendrian sequence
$$\rm{\pm stabilization \rightarrow isotopy \rightarrow \pm
destabilization}$$
corresponds to an elementary flype. Although the modeled isotopy in
Figure 8 of \cite{[BM4]}
is in the transversal setting, it is illustrative of the central issue.
The ``isotopy'' portion
of such a sequence can be highly complex---this particular modeled
isotopy requires repeated
uses of exchange moves. Thus, for the present note an argument is
needed to establish
that such a Legendrian sequence coming from moving between
differing classes having the same coordinate address on the
Legendrian mountain is realized by an elementary flype.\\


The outline of this note is as follows.  In \S\ref{section:standard contact structure} we state our main theorem, Theorem \ref{main theorem}, which establishes a Legendrian MTWS for the $(2,3)$ cabling of the $(2,3)$-torus knot.  In \S\ref{section:L+} we begin to work toward justifying this theorem by exhibiting a braided rectangular diagram of one of the Legendrian representatives of the $(2,3)$ cabling of the $(2,3)$-torus knot.  This requires us to develop a synthesis of convex tori, standard tilings and braided rectangular diagrams.  In \S\ref{section:proofoftheorem} we prove our main theorem.  We conclude with an appendix in \S 5 that further develops our synthesis of convex tori and tilings.



\noindent

\section{A knot-type specific Legendrian and transversal MTWS}
\label{section:standard contact structure}

\subsection{Background.} 
Consider $S^3$, viewed as the one-point compactification of $\mathbb{R}^3$.  In this context, the standard contact structure on $S^3$ can be thought of as the closure of the standard contact structure on $\mathbb{R}^3$, given in cylindrical coordinates as the kernel of the 1-form $dz+ r^2 d\theta$.  We denote this standard
tight contact structure by $\xi_{sym}$ in $\mathbb{R}^3$.
Given a topological knot type $\mathcal{K}$, we can restrict ourselves to look at representatives that are everywhere tangent to $\xi_{sym}$.  These are called Legendrian knots, and we say that two Legendrian knots are Legendrian isotopic if they can be connected by a 1-parameter family of Legendrian knots.  Similarly, we can restrict ourselves to look at representatives that are everywhere transverse to $\xi_{sym}$.  These are called transversal knots, and we say that two transversal knots are transversally isotopic if they can be connected by a 1-parameter family of transversal knots.\\
%

Both Legendrian and transversal knots have classical invariants, besides the topological knot type, that are preserved under Legendrian and transversal isotopies, respectively.  The Legendrian invariants are the rotation number, denoted by $r$, and the Thurston-Bennequin number, denoted by $tb$.  The transversal invariant is the self-linking number, denoted by $sl$.  A thorough discussion of these invariants, as well as general background to Legendrian and transversal knots, is provided in the excellent survey article by Etnyre found in \cite{[E1]}.\\

\subsection{Etnyre and Honda's Legendrian mountain range}
\begin{figure}[htbp]
	\centering
		\includegraphics[width=0.50\textwidth]{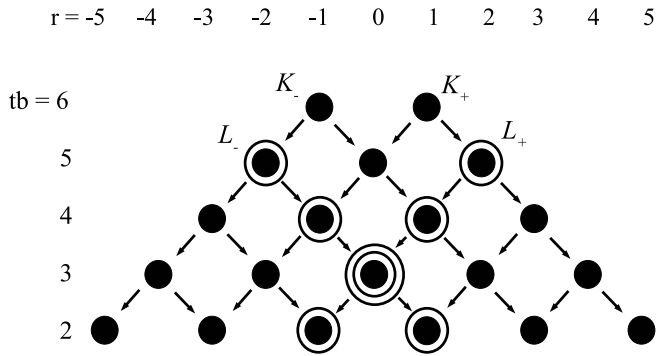}
	\caption{\small{The Legendrian mountain range for a $(2,3)$-cable of a $(2,3)$-torus knot.  A central dot and concentric circles represent multiple isotopy classes at a given value of $(r,tb)$.  Arrows down and to the left represent Legendrian negative stabilization; arrows down and to the right represent Legendrian positive stabilization.  $K_+$, $K_-$, $L_+$, and $L_-$ are defined in \cite{[EH]}.}}
	\label{fig:LegMountRangeb}
\end{figure}

For any topological knot type $\mathcal{K}$, one can represent the Legendrian isotopy classes as points on a two-dimensional grid, where the two coordinates are given by the values of $(r,tb)$ for that class.  For $\mathcal{K}$, there is a maximum value for the Thurston-Bennequin number, and thus this representation takes the shape of a {\em mountain range}; see Figure \ref{fig:LegMountRangeb}.  If there are multiple isotopy classes having the same value of $(r,tb)$, this can be represented by drawing circles around the central point, one for each multiple isotopy class.  Note that any mountain range is symmetric about the $r=0$ axis.  Arrows pointing down and to the left represent Legendrian negative stabilization, which we symbolically refer to as $S_-$; arrows pointing down and to the right represent Legendrian positive stabilization ($S_+$).\\

The Legendrian mountain range for a $(2,3)$-cable of a $(2,3)$-torus knot is shown in Figure \ref{fig:LegMountRangeb}, and was established by Etnyre and Honda in \cite{[EH]}.  The following structure is included in this mountain range:

\begin{itemize}

\item[1.]  At $tb=5$ and $r=2$, the outer circle represents $L_+$ while the inner dot represents $S_+(K_+)$; these are different Legendrian isotopy classes at the same values for the classical invariants.  A similar relationship holds for $L_-$ and $S_-(K_-)$.
\item[2.]  $S_+^k(L_-)$ is not Legendrian isotopic to $S_+^k(S_-(K_-))$ for any $k$; similarly $S_-^k(L_+)$ is not Legendrian isotopic to $S_-^k(S_+(K_+))$ for any $k$.  Also $S_+^2(L_-)$ is not Legendrian isotopic to $S_-^2(L_+)$.
\item[3.]  $S_-(S_+^k(L_-))=S_-(S_+^k(S_-(K_-)))$ and $S_+(S_-^k(L_+))=S_+(S_-^k(S_+(K_+)))$ for all $k$.

\end{itemize}

Of particular interest for this note is that in order to move from $L_+$ to the maximal $tb$ representative $K_+$, one must {\em first} stabilize, and {\em then} destabilize twice.  The main goal of this note is to show how to accomplish this movement toward maximal $tb$ without stabilization.


%
\subsection{The Legendrian MTWS for the $(2,3)$-cable of a $(2,3)$-torus knot}

We are now in a position to state our main theorem.

\begin{theorem}
\label{main theorem}
Let $\mathcal{K}_{(2,3)}$ be the $(2,3)$-cable of a $(2,3)$-torus knot.  Then Legendrian positive and negative destabilizations, along with elementary negative flypes, are sufficient to take a Legendrian representative of $\mathcal{K}_{(2,3)}$ to a representative at maximal Thurston-Bennequin number, modulo Legendrian isotopy.  In particular, we have the following:

\begin{itemize}
\item[1.]  $S_+^k(L_-)$ is related by an elementary negative flype to $S_+^k(S_-(K_-))$, for any $k$.
\item[2.]  $S_-^k(L_+)$ is related by an elementary negative flype to $S_-^k(S_+(K_+))$, for any $k$.

\end{itemize}
\end{theorem}

\begin{figure}[htbp]
	\centering
		\includegraphics[width=0.80\textwidth]{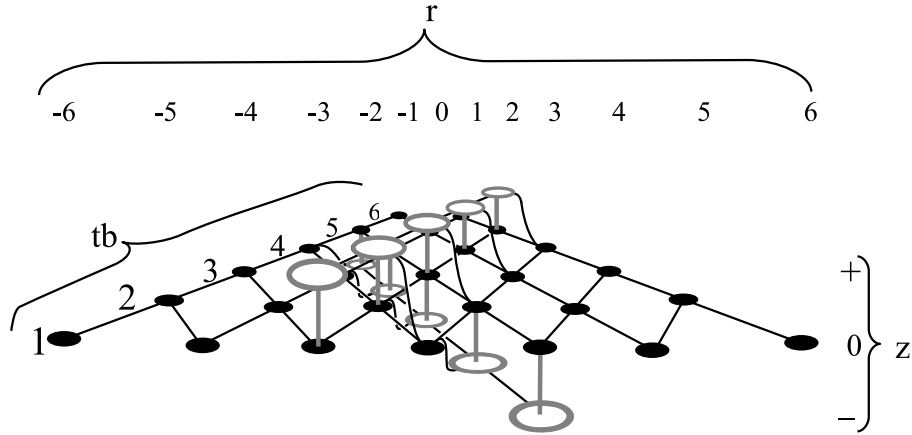}
	\caption{\small{Shown is the enhanced Legendrian mountain range for a $(2,3)$-cable of a $(2,3)$-torus knot.  The gray vertical lines indicate negative flypes performed on braided rectangular diagrams in two axes.}}
	\label{fig:EnhancedLegMtRngb}
\end{figure}

The structure of this theorem yields an enhanced Legendrian mountain range for the $(2,3)$-cable of a $(2,3)$-torus knot, depicted in Figure \ref{fig:EnhancedLegMtRngb}.  In this figure, black lines moving ``down'' the mountain range and to the right represent positive stabilization; similarly, black lines moving ``down'' the mountain range and to the left represent negative stabilization.  The new elements are the gray vertical lines connecting central dots to concentric circles either displaced above or below the central dots; these vertical lines represent elementary negative flypes.  The variable $z$ is a dummy variable; $z=0$ represents all stabilizations of $K_+$ or $K_-$, $z>0$ represents negative stabilizations of $L_+$, and $z<0$ represents positive stabilizations of $L_-$.

\subsection{The induced transversal MTWS for the $(2,3)$-cable of a $(2,3)$-torus knot}

The following theorem of Epstein, Fuchs, and Meyer connects the Legendrian classification of a knot type to its transversal classification.  $T_+$ denotes the positive transverse push-off of a Legendrian knot.

\begin{theorem}[Epstein, Fuchs, Meyer]
\label{theorem:push-off isotopy theorem}
Let $K_1$ and $K_2$ be two Legendrian knots.  Then the transversal knots $T_+(K_1)$ and $T_+(K_2)$ are transversally isotopic if and only if $S_-^m(K_1)$ and $S_-^n(K_2)$ are Legendrian isotopic for some $m$ and $n$ (where $m$ and $n$ could be zero).
\end{theorem}

This theorem allows us to obtain the transversal MTWS for the $(2,3)$-cable of a $(2,3)$-torus knot by taking positive transverse push-offs and collapsing the enhanced Legendrian mountain range to an {\em enhanced transversal trail}.  We thus obtain the following corollary:

\begin{corollary}
Let $\mathcal{K}_{(2,3)}$ be the $(2,3)$-cable of a $(2,3)$-torus knot.  Then negative braid destabilizations, along with elementary negative flypes, are sufficient to take a transversal representative of $\mathcal{K}_{(2,3)}$ to a representative at maximal self-linking number, modulo transversal isotopy.  In particular, $T_+(L_+)$ is related by an elementary negative flype to $T_+(S_+(K_+))$.
\end{corollary}

\begin{figure}[htbp]
	\centering
		\includegraphics[width=0.50\textwidth]{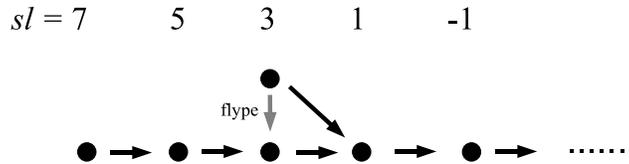}
	\caption{\small{Shown are the transverse isotopy classes for a $(2,3)$-cable of a $(2,3)$-torus knot.  The vertical arrow represents an elementary negative flype of braids, and the horizontal arrows represent negative braid stabilization.}}
	\label{fig:Transverseskylineb}
\end{figure}

Shown in Figure \ref{fig:Transverseskylineb} is the enhanced transversal trail depicting the structure in this corollary.



\textbf{Notation:}  Because we are interested in a particular knot type, for ease of notation the topological knot type of a $(2,3)$-torus knot will be denoted by $\mathcal{K}$, and the topological knot type of a $(2,3)$-cable of a $(2,3)$-torus knot will be denoted by $\mathcal{K}_{(2,3)}$.  The torus peripheral to $\mathcal{K}$ on which $\mathcal{K}_{(2,3)}$ resides will be denoted by $\mathcal{T}$.\\

\section{A braided rectangular diagram of $L_+$}
\label{section:L+}
We will prove Theorem \ref{main theorem} by explicitly demonstrating the necessary elementary negative flypes using braided rectangular diagrams for the Legendrian knots in question.  In order to do so, we must justify that the braided rectangular diagrams which we are using do actually correspond to representatives of the different Legendrian isotopy classes in Figure \ref{fig:LegMountRangeb}.  We begin that process in this section by constructing a braided rectangular diagram for $L_+$ in Figure \ref{fig:LegMountRangeb} and justifying that construction.\\

A note should be made here that the particular braided rectangular diagram for $L_+$ first appeared in \cite{[MM]}, along with the braided rectangular diagram resulting from a flype, as shown in Figure \ref{fig:MenascoFlypeb}.  These two diagrams were further studied in the context of knot Floer homology in \cite{[NOT]}, and from the combination of \cite{[MM], [EH], [NOT]} one can indirectly conclude that these two diagrams represent $L_+$ and $S_+(K_+)$, respectively.  However, a direct proof of why these diagrams represent $L_+$ and $S_+(K_+)$ has not been presented.  We do so in this section and the next to make explicit the structural connections between these independent lines of inquiry.

\subsection{$L_+$ as a Legendrian ruling}
Any torus in $(S^3,\xi_{sym})$ can be perturbed to be {\em convex}, meaning there exists a vector field everywhere transverse to the torus whose flow preserves the contact structure.  Recall that the characteristic foliation induced by the contact structure on a convex torus can be assumed to have a standard form, where there are $2n$ parallel {\em Legendrian divides} and a one-parameter family of {\em Legendrian rulings}. Parallel push-offs of the Legendrian divides gives a family of $2n$ {\em dividing curves}, referred to as $\Gamma$.  For a particular convex torus, the slope of components of $\Gamma$ is fixed; however, the Legendrian rulings can take on any slope other than that of the dividing curves by Giroux's Flexibility Theorem \cite{[G]}.\\

The knots in Figure \ref{fig:LegMountRangeb} are either Legendrian rulings or Legendrian divides on convex tori \cite{[EH]}.  For these convex tori, denoted by $\mathcal{T}$, two coordinate systems can be used.  One coordinate system, denoted by $\mathcal{C}_\mathcal{K}$, has a meridian of $\mathcal{T}$ having slope 0 and the preferred longitude of $\mathcal{T}$ having slope $\infty$.  The other coordinate system on $\mathcal{T}$, denoted by $\mathcal{C}_\mathcal{K}^\prime$, has a meridian having slope 0, while the curve having slope $\infty$ is found in the following manner:  Take the torus, peripheral to the unknot, on which $\mathcal{K}$ resides, and call it $T_0$.  Since $\mathcal{T}$ is peripheral to a representative of $\mathcal{K}$ on $T_0$, $\mathcal{T}$ will intersect $T_0$ in two parallel curves.  The slope of these curves on $\mathcal{T}$ is given the value $\infty$ in $\mathcal{C}_\mathcal{K}^\prime$.  As shown in \cite{[EH]}, $L_+$ is a Legendrian ruling on a convex torus that has two Legendrian divides of slope $-\frac{2}{11}$ in $\mathcal{C}_\mathcal{K}^\prime$.  $L_+$ intersects each of these Legendrian divides once in a positive intersection.  Moreover, the solid torus with boundary slope $-\frac{2}{11}$ is one which fails to thicken, meaning any solid torus containing it also has boundary slope $-\frac{2}{11}$.

\subsection{A convex torus representing $\mathcal{K}$ with $\textrm{slope}(\Gamma)=-\frac{2}{11}$}
Our goal in this subsection is to construct a solid torus representing $\mathcal{K}$ whose torus boundary has $\textrm{slope}(\Gamma)=-\frac{2}{11}$.  To do so, we will connect convex tori and Legendrian knots to the work of Menasco and Matsuda on standardly tiled tori and transversal knots represented as closed braids.\\

We first briefly review definitions found in \cite{[M1]} and \cite{[M3]}.  Consider $(\reals^3,\{z-axis\})\subset(S^3,\textbf{A})$, where $\textbf{A}$ is the axis of a transverse braid, $S^3=\reals^3 \cup \infty$, and $\textbf{A}=\{z-axis\} \cup \infty$.  In our case, this transverse braid will lie on a torus.  Let $(\rho,\theta,z)$ be the cylindrical coordinate system.  We denote the braid fibration by \textbf{H}=\{$H_\theta\mid0\leq\theta<2\pi$\}.  This will induce a singular braid foliation on the surface of the torus, where the singularities are either elliptic (where the torus intersects $\textbf{A}$) or hyperbolic (where the torus is tangent to a particular $H_\theta$).  Both the braid axis and the $H_\theta$'s have an orientation, and thus the singularities will be either positive or negative depending on whether the orientation of the torus agrees or disagrees with these orientations.  For our torus $\mathcal{T}$, singularities are joined by arcs which lie in a particular $H_\theta$, and each singularity will be connected to four other singularities via arcs.  The braid foliation on the torus is then said to be a {\em standard tiling}.  The particular standard tiling for the $\mathcal{T}$ that we will need is found in \cite{[M3]}, and is shown in Figure \ref{figure:2,3;1Knot}, along with a knot that is everywhere transverse to the foliation, and hence a braid.\\

\begin{figure}[htbp]
	\centering
		\includegraphics[width=0.60\textwidth]{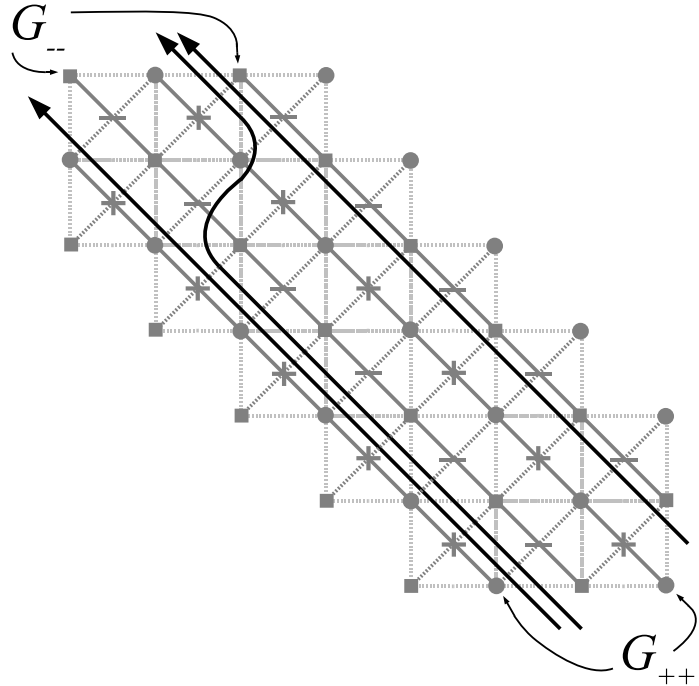}
	\caption{\small{Positive hyperbolic singularities are indicated by a $+$, negative hyperbolic singularities by a $-$.  Positive elliptic singularities are indicated by a gray dot, negative elliptic singularities by a gray square.  $K_{(2,3)}$ is indicated in black.}}
	\label{figure:2,3;1Knot}
\end{figure}

\begin{figure}[htbp]
	\centering
		\includegraphics[width=0.50\textwidth]{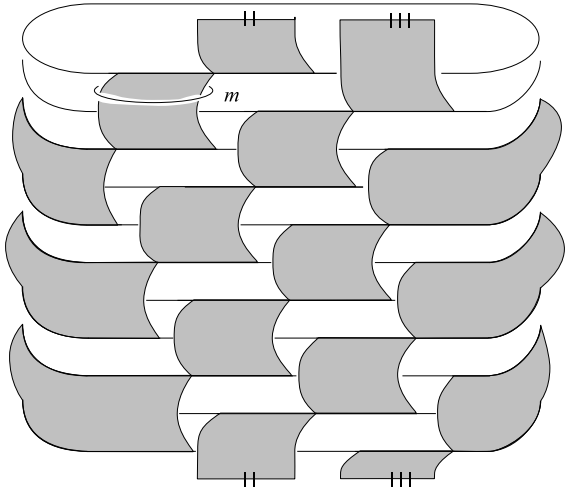}
	\caption{\small{A rectangular block presentation for the positive trefoil $K$, with blocks in gray and discs in white.  The boundary of a regular neighborhood $N$ of this collection of blocks and discs forms the torus $\mathcal{T}$.  A meridian curve of $\mathcal{T}$ is indicated by the letter $m$.}}
	\label{fig:PositiveTrefoilb}
\end{figure}

We want to see how $\mathcal{T}$ can be embedded in $S^3$.  We first construct the solid torus for which $\mathcal{T}$ is the boundary.  This solid torus can be represented using a rectangular block diagram, as described in \cite{[M3]}.  In particular, take a collection of discs of common radius whose centers are on the braid axis and which are parallel to the $xy$-plane.  We then attach to each disc a unique rectangular-shaped block whose bottom edge is on the boundary of the disc.  The top edges of the blocks are also attached to discs in a one-to-one fashion.  We do this so that the block-disc collection deformation retracts to a positive trefoil $K$.  A rectangular block diagram for $K$ is shown in Figure \ref{fig:PositiveTrefoilb}, where blocks are in gray and discs in white.  A regular neighborhood of the block-disc collection forms a solid torus $N$ whose boundary is $\mathcal{T}$.  Negative elliptic singularities are just below the centers of the discs, while positive elliptic singularities are just above the centers of the discs.    Negative hyperbolic singularities occur just to the left of the left edges of blocks, while positive hyperbolic singularities occur just to the right of the right edges of blocks.  A salient feature is that the left side of a block shares the same angular position as the right side of a block below it.  This implies that the associated negative hyperbolic singularity on $\mathcal{T}$ actually occurs before the positive hyperbolic singularity in the $\theta$-ordering.\\

\begin{figure}[htbp]
	\centering
		\includegraphics[width=0.50\textwidth]{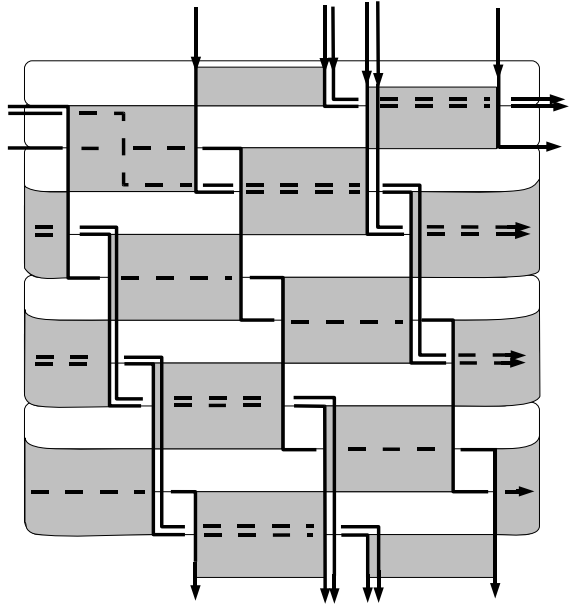}
	\caption{\small{$K_{(2,3)}$ on the rectangular block diagram of $\mathcal{T}$.  Vertical arcs going down the front of two consecutive blocks actually pass between the negative singularity that comes from the left edge of the top block, and the positive singularity that comes from the right edge of the bottom block.}}
	\label{fig:KNOTONPosTrefk=1b}
\end{figure}

Our knot of interest, $K_{(2,3)} \in \mathcal{K}_{(2,3)}$,  is a transverse braid on the surface of $\mathcal{T}$.  We can visualize this braid as being superimposed on the rectangular block diagram of $\mathcal{T}$.  This yields a braided rectangular diagram consisting of a collection of vertical and horizontal arcs, as defined in \cite{[M3]} and \cite{[MM]}.  The braided rectangular diagram for $K_{(2,3)}$ is shown in Figure \ref{fig:KNOTONPosTrefk=1b}.  Notice that $K_{(2,3)}$ has vertical arcs that go down along the front of the blocks, two blocks at a time, except for the one vertical arc in the upper left corner that passes behind one of the blocks.  The vertical arcs passing in front of two blocks at a time should be understood as running between the negative singularity that comes from the left edge of the top block, and the positive singularity that comes from the right edge of the bottom block.  It is clear that this knot has intersection number three with a meridian; by drawing the preferred longitude on the surface of $\mathcal{T}$, one can confirm that the knot has intersection number two with that longitude, and hence is a $(2,3)$-cabling.  This is the same knot pictured in Figure \ref{figure:2,3;1Knot}.\\

We now connect the braid foliation to the characteristic foliation induced by the contact structure.

\begin{lemma}
\label{lemma:braidfoliation equals charfoliation}
Suppose that the transverse braid $K_{(2,3)}$ lies on a standardly tiled torus $\mathcal{T}$, with the tiling induced by the braid fibration.  Then we may assume that $K_{(2,3)}$ lies on a standardly tiled torus $\mathcal{T}$, where the characteristic foliation is also a standard tiling.
\end{lemma}

\begin{proof}
In our block-disc collection, we choose a radius $r$ large enough so that at the boundary of the discs, the contact planes are $\epsilon$-close to being in the half planes $H_\theta$.  We then slightly tilt the blocks so that their sides are aligned with the contact planes at the large radius $r$.  Taking a neighborhood of this new block-disc presentation will give a solid torus on whose boundary the elliptic singularities will be positive-negative pairs where the $z$-axis intersects the discs, and the hyperbolic singularities will be positive-negative pairs that occur on the edges of each of the blocks. The characteristic foliation is thus a standard tiling.  We can do this while keeping the transverse isotopy class of the knot the same, and while maintaining the fact that the braid foliation is a standard tiling.
\end{proof}   

\begin{figure}[htbp]
	\centering
		\includegraphics[width=0.50\textwidth]{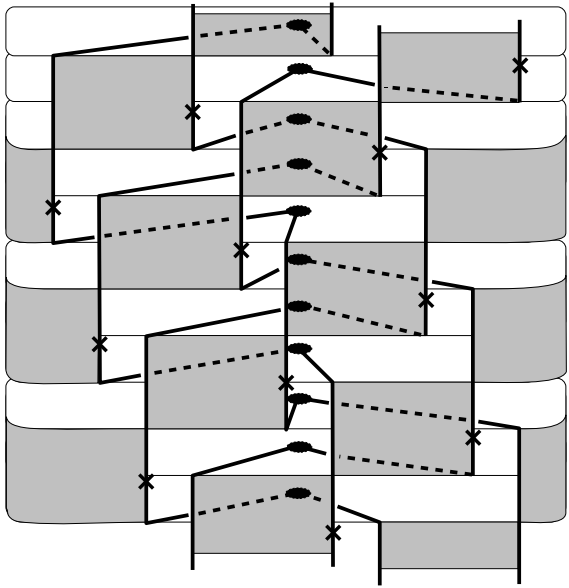}
	\caption{\small{$G_{++}$ superimposed on $\mathcal{T}$.  Elliptic singularities are dots; hyperbolic singularities are indicated by an x.}}
	\label{fig:G++onRectBlockk=1b}
\end{figure}

On a tiling, we can define four graphs, $G_{\epsilon\delta}$, that consist of elliptic singularities of parity $\epsilon$ connected by arcs that pass through hyperbolic singularities of parity $\delta$.  For a standard tiling, the components of $G_{++}$ and $G_{--}$ form a collection of an even number of parallel, homotopically non-trivial simple closed curves on the torus \cite{[M1]}.  Now $G_{++}$ and $G_{--}$ are piecewise Legendrian curves.  By a small isotopy of $\mathcal{T}$ near the braid axis, we may smooth out the corners and assume that $G_{++}$ and $G_{--}$ are Legendrian curves.  In Figure \ref{fig:G++onRectBlockk=1b}, we show $G_{++}$ superimposed on the rectangular block presentation for $\mathcal{T}$.  We have labelled the positive hyperbolic singularities with a black x, and the positive elliptic singularities with a black dot.\\

If we compare Figure \ref{fig:G++onRectBlockk=1b} with Figure \ref{fig:KNOTONPosTrefk=1b}, we can see that $K_{(2,3)}$ intersects $G_{++}$ right after the occurrence of the vertical arc that lies behind one of the blocks, and we can arrange things so this is the only intersection.  Similarly, if one imagines $G_{--}$ on $\mathcal{T}$, the only intersection of $K_{(2,3)}$ with $G_{--}$ occurs just before the occurrence of that same vertical arc.  This will be important in a coming subsection.  Now in the coordinate system $\mathcal{C}_\mathcal{K}^\prime$, $G_{++}$ intersects each meridian curve algebraically twice.  Moreover, the slope $\infty$ longitude intersects the top right corner of each block, and thus intersects $G_{++}$ once for each block.  The intersections of $G_{++}$ with $\infty$ are algebraically negative.  Since there are eleven blocks, the slope of $G_{++}$ in $\mathcal{C}_\mathcal{K}^\prime$ is $-\frac{2}{11}$.\\

The astute observer will notice that $\mathcal{T}$ is a convex torus with dividing curves that are parallel push-offs of $G_{++}$ and $G_{--}$, as such curves would separate the characteristic foliation of our torus into positive and negative regions.  We include the following proposition to formalize this:

\begin{lemma}
\label{prop:G++ and G-- the Legendrian divides}
Suppose $\mathcal{T}$ has a standard tiling that is the characteristic foliation.  Then we can isotop $\mathcal{T}$, rel $G_{++}$,$G_{--}$,$G_{+-}$, and $G_{-+}$ so that the resulting torus is standard form convex, with the components of $G_{++}$ and $G_{--}$ as the Legendrian divides and the components of $G_{+-}$ and $G_{-+}$ being Legendrian rulings.
\end{lemma}

\begin{proof}
The proof is an application of manipulation lemmas for the characteristic foliation found in \cite{[EF]}.  In short, $G_{++}$ and $G_{--}$ can be made to be simple closed curves of positive and negative singularities, respectively, by local manipulations of $\mathcal{T}$ that fix $G_{++}$ and $G_{--}$.  Moreover, this can be done while fixing $G_{+-}$ and $G_{--}$; the resulting torus thus has Legendrian rulings parallel to these two curves.
\end{proof}

As a consequence we have that the rectangular block presentation of $\mathcal{T}$ indeed represents a standard form convex torus with $\textrm{slope}(\Gamma)=-\frac{2}{11}$.  Moreover, there are two Legendrian divides. 

\subsection{The solid torus $N$ representing $\mathcal{K}$ fails to thicken}  From \cite{[EH]}, we know that the solid torus representing $\mathcal{K}$ with boundary slope $-\frac{2}{11}$ which fails to thicken has a particular complement in $S^3$.  Specifically, this complement consists of a regular neighborhood of a Legendrian representative of the Hopf link, with boundary consisting of two tori with boundary slopes $-\frac{1}{3}$ and $-\frac{1}{4}$, joined by a standard convex annulus with boundaries Legendrian representatives of $\mathcal{K}$.  In this subsection we show that the complement of a regular neighborhood of our block-disc collection is indeed this complement, and hence our block-disc collection represents a solid torus which fails to thicken.\\

\begin{figure}[htbp]
	\centering
		\includegraphics[width=0.50\textwidth]{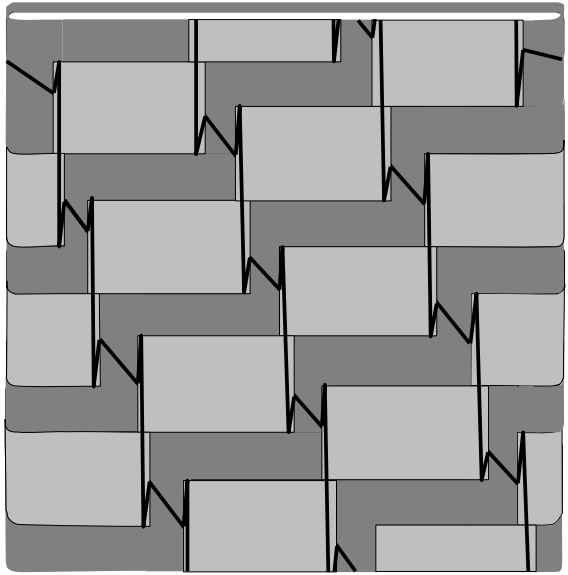}
	\caption{Shown in dark gray is the annulus $A$ connecting $N$ to itself; the blocks of $N$ are in light gray.  In black is one component of the dividing curves for the boundary of one of the solid tori bound by $N \cup N(A)$.  The boundary slope is $-\frac{1}{3}$.}
	\label{fig:Complement3BB}
\end{figure}

The proof is by picture.  We begin by letting $A$ be a standard convex annulus from $N$ to itself so that $N \cup N(A)$ is a thickened torus that bounds two solid tori representing unknots.  This is shown in Figure \ref{fig:Complement3BB},  where $A$ is shown in dark gray and the blocks from $N$ are in light gray.  Moreover, one component of the dividing curves for the boundary of one of the solid tori is shown in black, after edge rounding.  This is the solid torus that the reader sits inside; we are looking out of the solid torus and seeing its torus boundary.  From this perspective, it is evident that this boundary slope is $-\frac{1}{3}$.\\

\begin{figure}[htbp]
	\centering
		\includegraphics[width=0.50\textwidth]{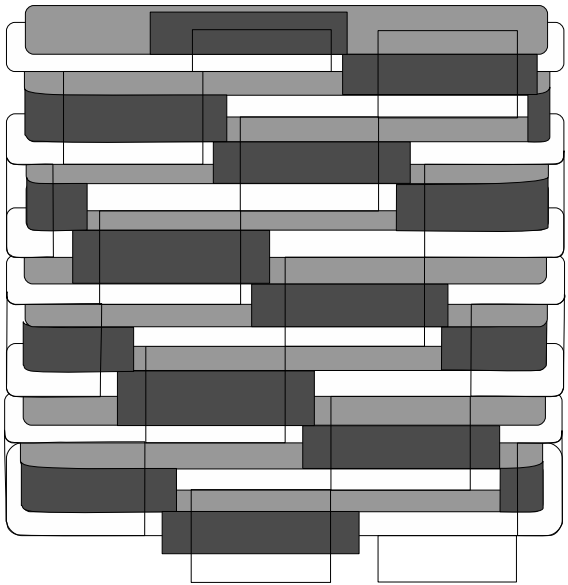}
	\caption{Shown is the second solid torus bound by $N \cup N(A)$ in blocks and discs with darker shades of gray.  The boundary slope of this solid torus is $-\frac{1}{4}$.}
	\label{fig:Complement2BB}
\end{figure}

Now note that the other solid torus bound by $N \cup N(A)$ has a block-disc representation as shown in Figure \ref{fig:Complement2BB}.  There, the blocks of $N$ are transparent, while the blocks and discs of the solid torus are two shades of gray.  It is evident that $G_{++}$, and hence the dividing curves, have slope $-\frac{1}{4}$.\\

This gives the correct solid tori in the complement of $N$; the standard convex annulus $\mathcal{A}'$ connecting them can be obtained by taking a longitudinal slice through the $I$-invariant neighborhood $N(A)$.  Thus $N$ is indeed the solid torus which fails to thicken from \cite{[EH]}.  We will refer to the braided rectangular diagram for the $(2,3)$-cabling of $K$, as shown in Figure \ref{fig:KNOTONPosTrefk=1b}, as $D$.\\ 

\textbf{Remark:}  There exists a solid torus representing $\mathcal{K}$ with boundary slope $-\frac{2}{11}$ that {\em does} thicken to a solid torus with boundary slope $-\frac{1}{5}$, i.e., one that is contained in a {\em standard neighborhood} of a $(2,3)$-torus knot at maximal Thurston-Bennequin value.  We include a block-disc presentation illustrating this in the Appendix.

\subsection{A braided rectangular diagram for $L_+$}
We now have the following lemma:

\begin{lemma}
The braided rectangular diagram $D$ represents the Legendrian isotopy class of $L_+$.
\end{lemma}

\begin{proof}
By Giroux's Flexibility Theorem, there is an isotopy of our convex torus such that the dividing curves remain fixed, but the resulting torus is still standard form convex, with Legendrian rulings that are $(2,3)$-cablings on our convex torus.  Moreover, we can accomplish this isotopy in the following way.  First imagine splitting our torus into two annuli bounded by the two Legendrian divides.  One of the annuli contains the portion of our knot that is parallel to the current Legendrian rulings (which recall are parallel translates of $G_{-+}$ and $G_{+-}$).  We leave this annulus fixed.  We isotop the other annulus using Giroux's Flexibility Theorem for surfaces with Legendrian boundary, so that the boundary of this annulus is still two curves of singularities, and the new rulings on the whole torus are $(2,3)$-torus knots.  In this way the rulings will intersect each dividing curve once, and each Legendrian divide once.  Because this isotopy kept the dividing curves fixed, we may still model our torus using a rectangular block presentation.  The portion of our knot that went behind the block will join up with a ruling that runs parallel to the Legendrian divides throughout the rest of its support.  Thus we would construct the braided rectangular diagram for our knot precisely the way we would construct the knot in \cite{[M3]}.  So the rectangular diagram in Figure \ref{fig:KNOTONPosTrefk=1b} from \cite{[M3]} is a rectangular diagram for the Legendrian representative of a $(2,3)$-cable of a $(2,3)$-torus knot that is a Legendrian ruling on a convex torus with Legendrian divides of slope $-\frac{2}{11}$, and which bounds a solid torus that fails to thicken.  Moreover, using formulas for the rotation number and Thurston-Bennequin number established in \cite{[MM]}, one can calculate that for the diagram in Figure \ref{fig:KNOTONPosTrefk=1b}, $r=2$ and $tb=5$.  Thus the diagram in Figure \ref{fig:KNOTONPosTrefk=1b} is a braided rectangular diagram for $L_+$ in \cite{[EH]}.
\end{proof}

\section{Proof of Theorem \ref{main theorem}}
\label{section:proofoftheorem}

\begin{proof}
We begin with $L_+$.  It is the outer circle at $(r,tb)=(2,5)$ in Figure \ref{fig:LegMountRangeb}, and we know from \cite{[EH]} that its Legendrian isotopy class does not Legendrian destabilize.  We claim that the inner dot at $(r,tb)=(2,5)$, which is $S_+(K_+)$, is represented by the braided rectangular diagram shown in \cite{[M3]} obtained by performing an elementary negative flype to the braided rectangular diagram of $L_+$.  This elementary negative flype can be thought of as forcing a Legendrian positive destabilization (which looks like a negative braid destabilization in the braided diagram) that can only be performed if a Legendrian positive stabilization occurs.  This results in $(r,tb)$ still being $(2,5)$.  $L_+$ and the diagram after the flype are shown in Figure \ref{fig:MenascoFlypeb}.\\

\begin{figure}[htbp]
	\centering
		\includegraphics[width=0.70\textwidth]{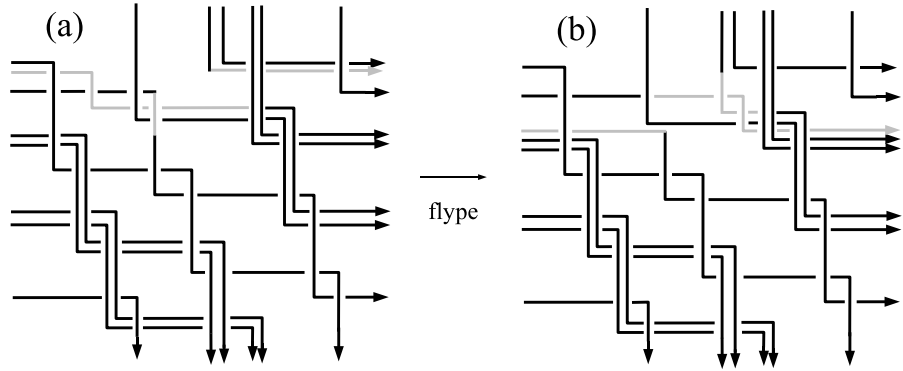}
	\caption{\small{$L_+$ is drawn in (a).  An elementary negative flype, seen as a Legendrian positive stabilization and Legendrian positive destabilization occurring on the gray arcs, is used to obtain the knot in (b).}}
	\label{fig:MenascoFlypeb}
\end{figure}

We now need to show that after the flype, the resulting knot destabilizes after a sequence of Legendrian isotopies, thus proving that it is indeed $S_+(K_+)$.  In particular, consider Figure \ref{fig:LegPosStabb}.  Part (a) is the knot obtained from $L_+$ following the elementary negative flype.  Indicated in (a) is a vertical arc in gray toward the left of the diagram.  If we move this to the right via a Legendrian isotopy, we can perform a Legendrian flip on the gray horizontal arc and then slide this horizontal arc upward over a black horizontal arc.  After reversing the original flip, we obtain the knot in (b).  In (b), we now focus in on the dashed arcs.  We can slide the right dashed vertical arc to the left, and then do a Legendrian flip on the horizontal dashed arc.  We can then slide that horizontal arc down over two black horizontal arcs.  After reversing the flip, we obtain (c).  In (c), there is now a Legendrian positive destabilization indicated by the shaded box.  If we perform this destabilization, we arrive at $K_+$ in the Legendrian mountain range in Figure \ref{fig:LegMountRangeb}.   Taking (c) in Figure \ref{fig:LegPosStabb} and doing a Legendrian flip and some Legendrian isotopies, we obtain $K_+$ as drawn in Figure \ref{fig:K+b}.\\

\begin{figure}[htbp]
	\centering
		\includegraphics[width=0.70\textwidth]{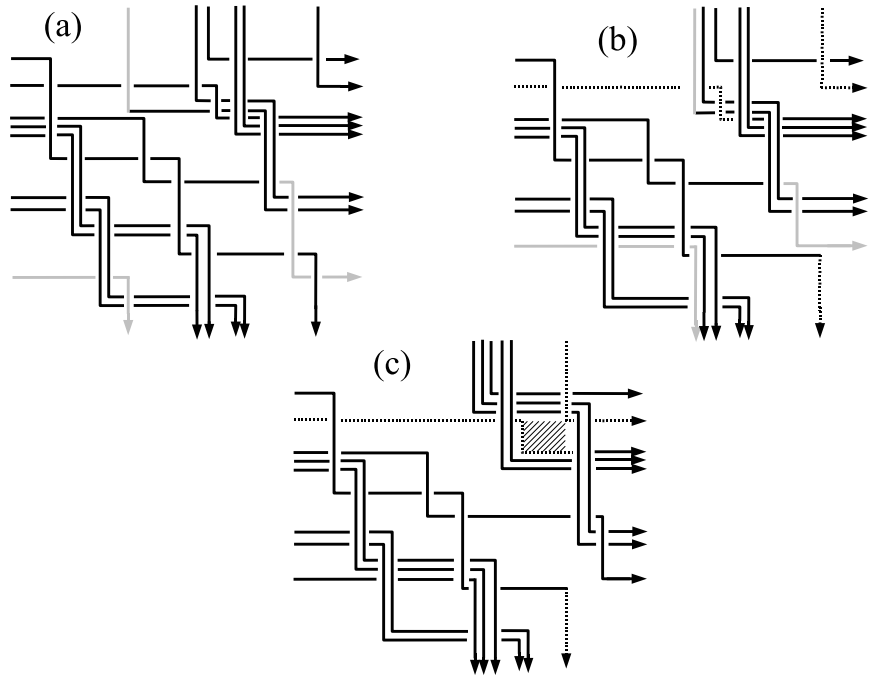}
	\caption{\small{Moving from (a) to (b) to (c) can be accomplished by a sequence of Legendrian moves that reveal a Legendrian positive destabilization, indicated by the shaded box.}}
	\label{fig:LegPosStabb}
\end{figure}

\begin{figure}[htbp]
	\centering
		\includegraphics[width=0.40\textwidth]{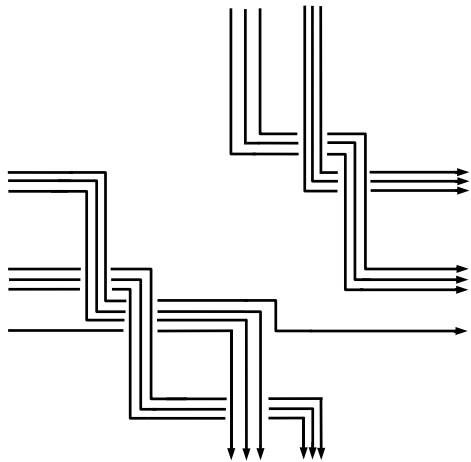}
	\caption{\small{Shown is $K_+$.}}
	\label{fig:K+b}
\end{figure}

This proves that $L_+$ and $S_+(K_+)$ are related by an elementary negative flype.  Moreover, if we perform $k$ Legendrian negative stabilizations on $L_+$ and $S_+(K_+)$ away from the support of the flype, we can perform them so that the rectangular diagrams of $S_-^k(L_+)$ and $S_-^k(S_+(K_+))$ differ by an elementary negative flype.\\

To complete the proof of the theorem we turn our attention to $L_-$ and $S_-(K_-)$.  To obtain braided rectangular diagrams for these knots we examine a general braided rectangular diagram $K$ at $(r,tb)$.  We take this braided rectangular diagram of $K$, and imagine it as projected onto a square.  We then flip the square along the diagonal that runs from the top right to the bottom left.  This is a topological isotopy of $K$.  One then reverses the orientation, yielding vertical arcs that pass over horizontal arcs, with vertical arcs pointing down and horizontal arcs pointing to the right.  This gives the braided rectangular diagram for the knot at $(-r,tb)$.\\
  
\begin{figure}[htbp]
	\centering
		\includegraphics[width=0.65\textwidth]{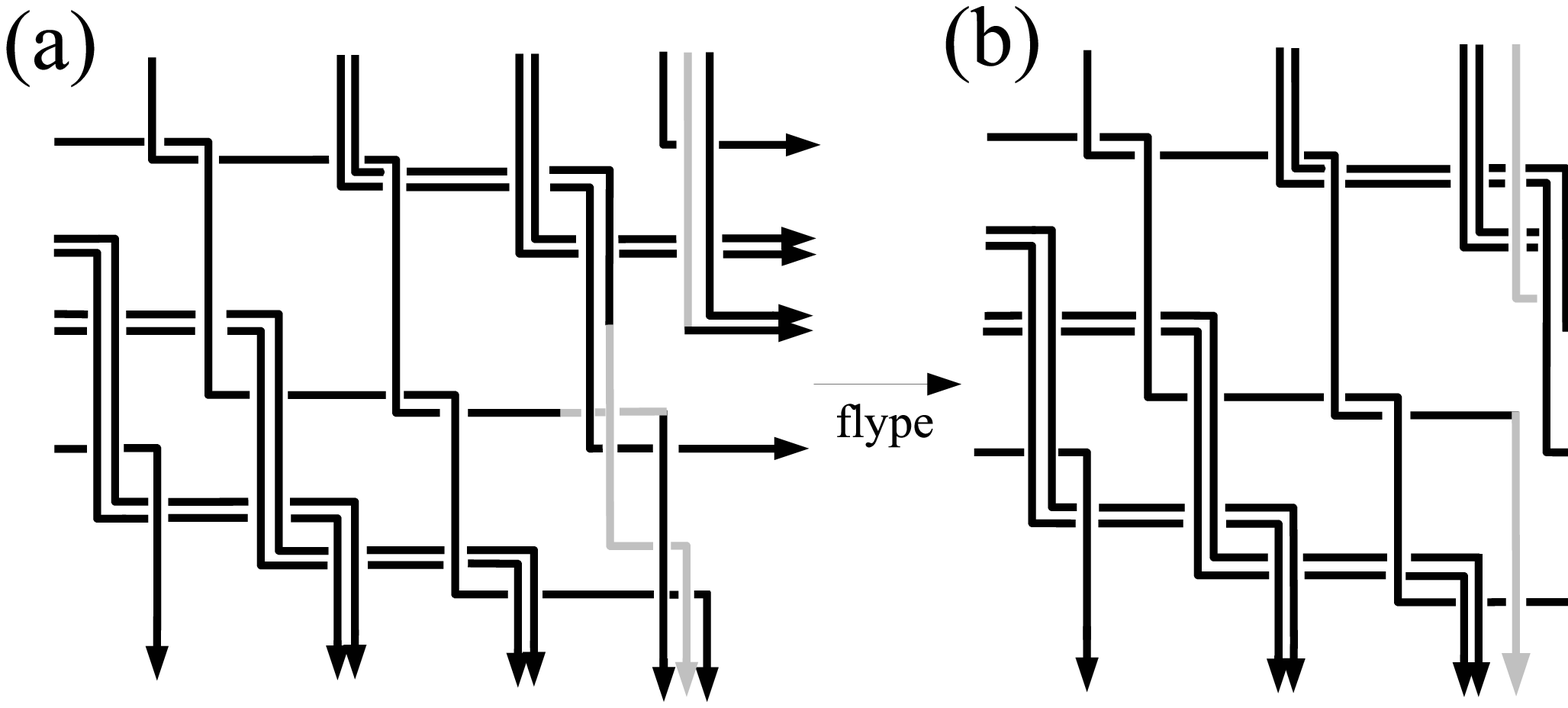}
	\caption{\small{In (a) is $L_-$; in (b) is $S_-(K_-)$.  They are related by an elementary negative flype in the horizontal braid axis.}}
	\label{fig:L-andflypeb}
\end{figure}

In particular, we show $L_-$ and $S_-(K_-)$ in Figure \ref{fig:L-andflypeb}.  To get from $L_-$ to $S_-(K_-)$, one forces a Legendrian negative destabilization which results in an accompanying Legendrian negative stabilization, as indicated in the gray arcs.  The move from $L_-$ to $S_-(K_-)$ is an elementary negative flype in the horizontal braid axis.  A flip of the vertical arcs is a Legendrian isotopy, and so one can show as before that the knot in (b) Legendrian negatively destabilizes, and in fact is $S_-(K_-)$.  Moreover, we again have that $S_+^k(L_-)$ and $S_+^k(S_-(K_-))$ are related via an elementary negative flype by performing the positive stabilizations outside the support of the original flype.  This concludes the proof.
\end{proof}

\section{Appendix -- A convex torus that admits a thickening}

The goal of this section is to provide a block-disc presentation for a solid torus representing the $(2,3)$-torus knot with boundary slope $-\frac{2}{11}$ that {\em thickens} to a standard neighborhood of a $(2,3)$ torus knot at maximal Thurston-Bennequin value.  To do this, we will (1) begin with our previous solid torus of slope $-\frac{2}{11}$ that {\em fails} to thicken; (2) show how to perturb it (without thickening) to obtain a solid torus with slope $-\frac{1}{5}$; and (3) show how to {\em thin} the solid torus with slope $-\frac{1}{5}$ to a new one with slope $-\frac{2}{11}$.\\

\begin{figure}[htbp]
	\centering
		\includegraphics[width=0.70\textwidth]{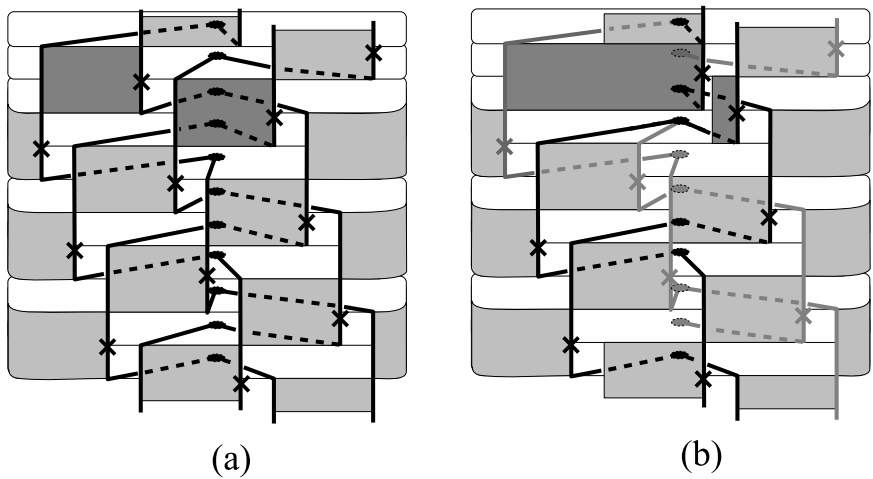}
	\caption{In (a) is the solid torus representing $\mathcal{K}$ with boundary slope $-\frac{2}{11}$ that fails to thicken.  A Legendrian divide is indicated in black.  In (b) is the solid torus representing $\mathcal{K}$ with boundary slope $-\frac{1}{5}$.  A Legendrian divide is indicated in black.  Elliptic and hyperbolic singularities removed by Giroux elimination are in darker gray.}
	\label{fig:211to15BB}
\end{figure}

In Figure \ref{fig:211to15BB} we show in (a) the solid torus that fails to thicken, as established in \S 3.  The graph $G_{++}$, which recall represents a Legendrian divide, is indicated in black, connecting elliptic and hyperbolic singularities.  The slope of this curve is $-\frac{2}{11}$ in the $\mathcal{C}'$ framing.  Note that two blocks are in a darker shade of gray.  We then make one of these dark gray blocks thinner, and one thicker, in order to arrive at the block-disc presentation of a solid torus in (b).  There the elliptic and hyperbolic singularities of $G_{++}$ are indicated, but note that $G_{++}$ forms a simple closed curve with two trees extending off of it; the simple closed curve is indicated in black, the trees in gray.  We can perturb this convex torus to be in standard form by using Giroux elimination along the gray trees to remove singularities, resulting in a black Legendrian divide with slope $-\frac{1}{5}$ in the $\mathcal{C}'$ framing.\\

\begin{figure}[htbp]
	\centering
		\includegraphics[width=0.70\textwidth]{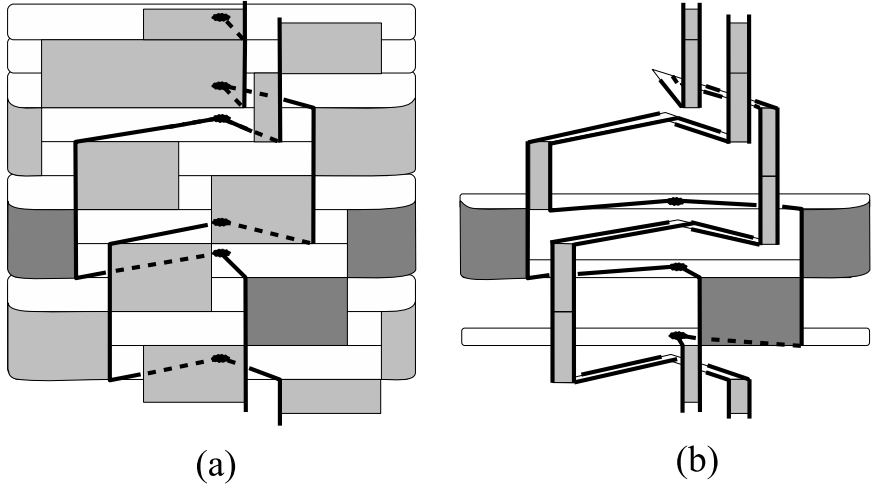}
	\caption{In (a) is the solid torus with boundary slope $-\frac{1}{5}$; a Legendrian divide is in black, and two blocks are in dark gray.  To get to (b), we fix these two blocks and their adjacent discs, and then thin the other blocks and discs, totally removing some discs.  The result is a standard convex torus with slope $-\frac{2}{11}$ that admits a thickening.}
	\label{fig:thinning15BB}
\end{figure}

We now show how to {\em thin} this solid torus to obtain one with boundary slope $-\frac{2}{11}$.  Refer to Figure \ref{fig:thinning15BB}, where in (a) we have redrawn the solid torus with boundary slope $-\frac{1}{5}$; a Legendrian divide is in black, and two blocks are in dark gray.  To get to (b) we fix these two blocks and their adjacent discs, and then {\em thin} the other blocks and discs, eliminating some discs entirely.  Then with a small perturbation of the solid torus, which we can accomplish by thinning, we can obtain the Legendrian divide indicated in black in (b), which has slope $-\frac{2}{11}$.  Reversing this thinning process then results in a thickening of this solid torus with boundary slope $-\frac{2}{11}$ to a standard neighborhood of a $(2,3)$-torus knot at maximal Thurston-Bennequin value.\\

\bigskip
\begin{itemize}
     \item[ ]\small{\scshape{University at Buffalo, Buffalo, NY}
     \item[ ]{\em E-mail addresses}:  \ttfamily{djl2@buffalo.edu}, \ttfamily{menasco@buffalo.edu}}
\end{itemize}
\bigskip

\end{document}